\documentclass[12pt]{article}
\usepackage[utf8]{inputenc}
\usepackage[english]{babel}
\usepackage{graphicx}
\usepackage{vmumath}
\usepackage{amsmath}
\usepackage{amsthm}
\usepackage{avo_book}
\usepackage{tikz}
\usepackage{caption}
\usepackage{booktabs}
\usepackage{subcaption}
\usepackage[a4paper,left=15mm,right=15mm,top=30mm,bottom=20mm]{geometry}
\noaftermath
\begin{document}

\bibliographystyle{unsrt}

\title{Enumeration of 4-regular one-face maps}

\author{
Evgeniy Krasko \qquad  Alexander Omelchenko\\
\small St. Petersburg Academic University\\
\small 8/3 Khlopina Street, St. Petersburg, 194021, Russia\\
\small\tt \{krasko.evgeniy, avo.travel\}@gmail.com
}

\begin{abstract}
We give explicit formulas enumerating 4-regular labelled and unlabelled one-face maps.

\end{abstract}

\maketitle

\section{Introduction}

By an one-face (or unicellular) map $\M$ on a compact connected orientable surface which is fully characterized by its genus $g$, we will denote an imbedding of a connected graph $G$, loops and multiple edges allowed, into a compact and oriented surface $X$, such that $G$ is as a subset of $X$ and complement $X\setminus G$ is homeomorphic to a topological polygon. This complement is the only face of the map $\M$, having $v$ vertices (points on the surface $X$) and $n$ edges (nonintersecting curves on the surface that have no common points other than the vertices of the graph). 

$4$-regular maps play an important role in different fields of mathematics. They have been used for the knot problem in low-dimensional topology \cite{LandoZvonkin_engl}, for rectilinear embedding in VLSI \cite{Liu}, for the Gaussian Crossing Problem \cite{Liu}, for the folding of RNA interaction structures \cite{RNA}, and for enumerating some other kinds of maps. In its turn, one-face maps are extensively used in a number of instances in pure mathematics including finite type invariants of knots and links \cite{Bar-Natan},\cite{Kontsevich}, the representation theory of Lie algebras \cite{Campoamor}, the geometry of moduli spaces of flat connections on surfaces \cite{Andersen}, mapping class groups \cite{Andersen2} and in applied mathematics including codifying the pairings among nucleotides in RNA molecules \cite{Reidys_book} and data structures analysis \cite{Flajolet}. The present paper is devoted to the problem of enumerating 4-regular one-face maps or, equivalently, objects dual to these maps, namely $1$-vertex quadrangulations of a surface of genus $g$. 

The first work on the subject of enumerating rooted maps on genus $g$ surfaces was the work of Walsh and Lehman \cite{Walsh_Lehman}. Using Tutte's approach for enumerating planar maps \cite{Tutte_Census}, \cite{Tutte_Enum}, the authors derived the recurrence relation for the numbers of rooted maps and calculated the first terms of the corresponding sequnces. In addition, that paper contained an explicit expression for the number $\epsilon_g(n)$ of one-face maps with $n$ edges on a surface of genus $g$, as well as a formula for the number of unicellular maps of genus $g$ with prescribed vertex degrees. However, the combinatorial sense of these formulas remained unclear. The first success in deriving similar results combinatorially was the work \cite{Chapuy}. Chapuy used an original approach based on the reduction of an one-face map to an one-face map of a lower genus. He obtained a new recurrence relation for the numbers $\epsilon_g(n)$ and gave an elegant combinatorial interpretation of it. In the same paper it was shown how to use this technique to enumerate special kinds of maps, for example, cubic one-faced maps. In the first part of the present work we apply this approach to enumerate $4$-regular rooted one-face maps of genus $g$. 

In the eighties of the last century the first papers \cite{Liskovets_85}, \cite{Wormald} devoted to enumerating unrooted planar maps appeared. Ideas of the work \cite{Liskovets_85} were significantly developed and improved by Mednykh and Nedela in \cite{Mednykh_Nedela}. Mednykh and Nedela introduced a concept of a map on an orbifold, i.e. a map on a quotient of a surface under a finite group of automorphisms. They enumerated unrooted maps on an orientable surface of a given genus $g$ with a given number of edges $n$. To accomplish this, they derived a formula for the number of order-preserving epimorphisms from the fundamental group of an orbifold onto the cyclic group $Z_l$. In the second part of the present work we use the approach described in \cite{Mednykh_Nedela} and the results obtained in the first part to enumerate $4$-regular unrooted unicellular maps of a given genus $g$.  

\section{Statement of the Problem}

\begin{figure}[ht]
\centering
    	\includegraphics[scale=1.3]{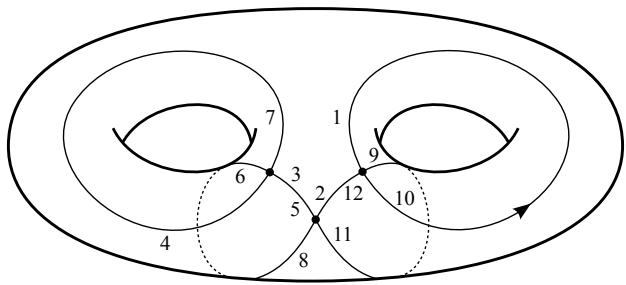}
	\caption{An one-face 4-regular map}
\label{fig:intertwined}
\end{figure}

For the problem of enumeration it is convenient to use the combinatorial definition of an unicellular map. Namely, any one-face map $\M$ having $n$ edges and $k$ vertices could be specified by the triple $(H,\alpha,\sigma)$, where $H$ describes the set of semi-edges of $\M$, $|H|=2n$. The vertices, edges and the only face of $\M$ are defined by the cycles of the permutations $\alpha$, $\sigma$ and $\gamma=\alpha\sigma$, where $\alpha$ is a fix-point-free involution, $\gamma$ is a cyclic permutation of length $2n$, and $\sigma$ is a collection of $k$ cycles $\omega_i$ corresponding to the vertices of $\M$. For example, figure \ref{fig:intertwined} shows an one-face map of genus $2$ with $n=6$ edges and $k=3$ vertices. This map is described by the permutations
$$
\alpha=(1,10)(2,12)(3,5)(4,7)(6,8)(9,11),\qquad \sigma=(1,12,10,9)(2,5,8,11)(3,7,6,4),
$$
$$
\gamma=\alpha\cdot \sigma=(1,2,3,4,5,6,7,8,9,10,11,12). 
$$
Euler’s characteristic formula states that the number $k$ of vertices and the number $n$ of edges uniquely determine the genus $g$ of the surface into which a graph is embedded:
\begin{equation}
\label{eq:Euler}
n-k=2g-1.
\end{equation}
This formula makes it possible to determine for given $\alpha$ and $\sigma$ the genus $g$ of an arbitrary map $\M$ and for given $g$ and $n$ compute the number of vertices $k$.

From the handshaking lemma it follows that for a $4$-regular map the number of vertices is $v=n/2$. Taking into account the formula (\ref{eq:Euler}), we obtain that the number $n$ of edges and the number $k$ of vertices are given by
\begin{equation}
\label{eq:4-regular}
n=4g-2,\qquad \qquad k=2g-1.
\end{equation}

A map is rooted if a semi-edge or a dart (an edge-vertex incidence pair) is distinguished as its root. By counting maps we mean counting equivalence classes of maps under orientation-preserving homeomorphism. For rooted maps any homeomorphism should preserve the distinguished oriented edge. In this case a homeomorphism preserves all the darts \cite{Tutte_Enum}, so that rooted maps can be counted without considering symmetries.  That is why this paper starts with the enumeration of rooted $4$-regular unicellular maps. Unrooted maps are enumerated in the second part. We also include a table listing the numbers of rooted and unrooted $4$-regular unicellular maps counted by genus.

\section{Enumeration of rooted $4$-regular one-face maps} 

Following \cite{Chapuy}, we denote by $r$ the root edge of a map $\mathfrak{m}=(H,\alpha,\sigma)$ and introduce a linear order on the set $H$ of its semi-edges as follows:
$$
r<\gamma(r)<\gamma^2(r)<\ldots<\gamma^{2n-1}(r).
$$
In other words, the ordering of $H$ is induced by the traversal order of semi-edges along the only face of $\mathfrak{m}$ performed with edges on the left (figure \ref{fig:intertwined}). This ordering, in turn, implies the ordering of semi-edges in a cycle
$$
\omega_i=(h^1,h^2,\ldots,h^m)
$$
of a permutation $\sigma$ that describes some vertex $v_i$ of a map $\mathfrak{m}$. An important observation made in \cite{Chapuy} is the following. Consider three arbitrary semi-edges $a_1$, $a_2$ and $a_3$ that belong to a common vertex $v_i$. In the planar case $(g = 0)$ when traversing along a planar tree $\mathfrak{m}$ performed with edges on the left, these semi-edges in the permutation $\sigma$ are ordered in the same way as in the permutation $\gamma$. In the case $g > 0$ necessarily exists a vertex $v_i$ and a triple $(a_1, a_2, a_3)$ of semi-edges belonging to it, such that their order in $\omega_i$ (for example, $(a_1, a_2, a_3)$) is opposite from that in the cycle of $\gamma$ (namely, $(a_1, a_3, a_2)$). We say that such a triple of semi-edges is intertwined in $v_i$. Using these semi-edges, one can cut the vertex into three new vertices, getting some new one-face map $\mathfrak{m'}=(H,\sigma',\alpha')$ on a surface of genus $g-1$ with the same number of edges. Conversely, for a map $\mathfrak{m'}$ of genus $g-1$ one can choose three of its arbitrary semi-edges $a_1<a_2<a_3$ belonging to three different vertices and glue these vertices into one new vertex $v$. As a result of this gluing operation a new one-face map $\mathfrak{m} $ of genus $g$ is obtained and its semi-edges $a_1$, $a_2$ and $a_3$ are intertwined at $v$. In \cite{Chapuy} it is shown that these two operations are inverses of each other. 

Moreover, a stronger result is proved in \cite{Chapuy}.

\begin{defin}
A semi-edge $h \in H$ is said to be a trisection if $\sigma(h)<h$, but $\sigma(h) \neq \min v(h)$, $v(h)$ being the vertex incident to $h$, $\min v(h)$ being the minimal semi-edge incident to $v(h)$.
\end{defin}

\begin{lemm}[Chapuy]
Every one-face map $\mathfrak{m} = (H, \alpha, \sigma) $ of genus $g$ has exactly $2g$ trisections.
\end{lemm}

Presence of a trisection in a vertex is equivalent to presence of some intertwined triple in it. The problem is that more intertwined triples than trisections could exist. However, Chapuy has showed a way of mapping each trisection to a unique intertwined triple. Moreover, $2g$ trisections of a map $\mathfrak{m}$ are divided into two groups depending on the result of cutting the vertex at the corresponding intertwined triple. The first one (so called trisection of type I) leads to a map of genus $g-1$ with three distinguished vertices. The second one (so called trisection of type II) corresponds to a map of genus $g-1$ with two distinguished vertices $v_1$, $v_2$ and a vertex $v(\tau)$ with a distinguished trisection $\tau$, such that $\min(v_1)<\min(v_2)<\min(v(\tau))$. These vertex-cutting operations are bijections between the corresponding classes of maps. 

We apply the approach described above for the enumeration of $4$-regular rooted one-face maps. Consider an one-face $(1 \div 4) $ map of genus $g$, i.e. some map $\mathfrak{m}$ having only the vertices of degrees one or four. Let $n$ be the number of edges in the map $\mathfrak{m}$, $s$ be the number of leafs, and $k$ be the number of vertices of degree $4$. From the Euler's formula and the handshaking lemma it follows that these numbers are connected with the genus $g$ of the map $\mathfrak{m}$ by the relations
\begin{equation}
n = 3k + 1-2g, \qquad\qquad s = 2k + 2-4g.
\label{n_s}
\end{equation}
Let $\epsilon_g^{(1 \div 4)}(k)$ be the number of one-face $(1 \div 4)$-maps of genus $g$ with $k$ vertices of degree $4$. For each of these maps there are $2g$ ways to distinguish a trisection, so the total number of rooted $(1 \div 4)$-maps with a distinguish trisection is equal to $2g \cdot \epsilon_g^{(1 \div 4)} (k)$.

Note that for $(1 \div 4)$ maps any trisection $\tau$ necessary belongs to type I. Indeed, a trisection of the other type can only appear when at least one of three vertices obtained by cutting the vertex with the trisection $\tau$ has a degree greater than or equal to three. So, cutting any such map $\mathfrak{m}$ with a distinguished trisection yields a rooted one-face $(1 \div 2 \div 4)$-map $\mathfrak{m'}$ with three distinguished vertices having degrees $1$, $1$ and $2$. In this map we will have $k-1$ vertices of degree $4$, the only vertex of degree $2$ and $s + 2$ vertices of degree $1$. One could erase a vertex of degree $2$, replacing two of its incident edges by a single edge to get a $(1 \div 4)$-map. However, this wouldn't establish a bijection because one of the erased semi-edges incident to a vertex of degree 2 might have been a root.

To obtain a counting formula, consider an arbitrary $(1 \div 4)$-map of genus $g-1$ with $(n-1)$ edges, $(k-1)$ vertices of degree $4$ and $(s + 2)$ leafs. For this map there are $\BCf{s + 2}{2}$ ways to distinguish two leafs and $(n-1)$ ways to distinguish an edge. If we subdivide a distinguished edge into two edges by putting a vertex of degree $2$ in its middle, two semi-edges incident to the new vertex are guaranteed not to be a root. To remedy this situation we use the principle of double counting. Namely, one could count maps having two roots, say, red and green, with a red root not be incident to a vertex of degree 2, in two ways: by choosing a red root in a green-rooted map or by choosing a green root in a red-rooted map. This results in an additional factor ${2n}/{(2n-2)}$ in the following number of ways to get a rooted $(1 \div 4)$-map of genus $g$ with $n$ edges, $k$ vertices of degree $4$ and $s$ leafs from a rooted $(1 \div 4)$ map of genus $g-1$ with $(n-1)$ edges, $(k-1)$ vertices of degree $4$ and $(s + 2)$ leafs:
$$
\BCf{s+2}{s}\cdot (n-1)\cdot\dfrac{2n}{2n-2}=n\cdot \BCf{s+2}{s}.
$$

The described correspondence between the maps of genus $g$ and genus $g-1$ enables us to write the following recurrence relation for the number $\epsilon_g^{(1\div 4)}(k)$ of one-face rooted $(1\div 4)$-maps with $n=3k+1-2g$ edges:
$$
2g\cdot \epsilon_g^{(1\div 4)}(k)=n\cdot \BCf{s+2}{s}\cdot \epsilon_{g-1}^{(1\div 4)}(k-1)\qquad \Longrightarrow\qquad
 \epsilon_g^{(1\div 4)}(k)=\dfrac{n}{2g}\cdot \BCf{s+2}{s}\cdot \epsilon_{g-1}^{(1\div 4)}(k-1).
$$

This recurrence could be easily unwinded to the genus $g=0$:
$$
\epsilon_g^{(1\div 4)}(k)=\dfrac{1}{2^g\,g!}\cdot\dfrac{n!}{(n-g)!}\cdot \BCf{s+2}{s}\cdot \BCf{s+4}{s+2}\cdot\ldots\cdot
\BCf{s+2g}{s+2g-2} \cdot\epsilon_{0}^{(1\div 4)}(k-g)\qquad \Longleftrightarrow
$$
$$
\Longleftrightarrow\qquad  
\epsilon_g^{(1\div 4)}(k)=\dfrac{1}{2^g\,g!}\cdot\dfrac{n!}{(n-g)!}\cdot \dfrac{(s+2g)!}{s!\,2^g} \cdot\epsilon_{0}^{(1\div 4)}(k-g).
$$

It remains to obtain an expression for the numbers $\epsilon_0^{(1\div 4)}(k_t)$ of rooted $(1\div 4)$-trees on a sphere that have $n_t$ edges, $k_t$ vertices of degree $4$ and $s_t$ leafs. From the handshaking lemma and the Euler's equation we get
$$
n_t=3k_t+1,\qquad\qquad s_t=2k_t+2.
$$
It's known that the numbers of $4$-trees having $k_t$ inner vertices and rooted at one of the leafs are Fuss-Catalan numbers
$$
C_4(k_t)=\dfrac{1}{k_t}\BCf{3k_t}{k_t-1}=\dfrac{1}{k_t}\BCf{n_t-1}{k_t-1}=\dfrac{(n_t-1)!}{k_t!(s_t-1)!}.
$$
Multiplying $C_4(k_t)$ by the number $2n_t$ of semi-edges, we get the number of $(1\div 4)$-trees having two roots. One of these roots is a leaf, and the other is a semi-edge. The same number could be obtained by multiplying $\epsilon_0^{(1\div 4)}(n_t)$ by the number $s_t$ of leafs. This double-counting arguments yields  
$$
\epsilon_0^{(1\div 4)}(k_t)=\dfrac{2\,n_t!}{k_t!\,s_t!}.
$$
Note that this expression is valid for any $(1\div d)$-trees. Taking into account that $n_t=n-g$, $k_t=k-g$ and $s_t=s+2g$, we conclude that the number of  $(1\div 4)$-maps on the surface of genus $g$ is given by
\begin{equation}
\label{eq:num_1_4_maps}
\epsilon_g^{(1\div 4)}(k)=\dfrac{n!\,(s+2g)!}{4^g\,g!\,s!\,(n-g)!}\cdot \dfrac{2\,(n-g)!}{(k-g)!\,(s+2g)!}=
\dfrac{2\,n!}{4^g\,g!\,s!\,(k-g)!},
\end{equation}
where $n$ and $s$ are related to $k$ and $g$ by (\ref{n_s}). Substituting $s=0$, $k=2g-1$ and $n=4g-2$ into (\ref{eq:num_1_4_maps}) we get the formula for the numbers $\epsilon_g^{(4)}$ of rooted $4$-regular one-face maps
\begin{equation}
\label{eq:num_4_maps}
\epsilon_g^{(4)}=\dfrac{2\,(4g-2)!}{4^g\,g!\,(g-1)!}; \qquad \text{sequence} \quad 1,45,9450,4729725,\ldots
\end{equation}

\section{Enumeration of unrooted $4$-regular one-face maps}

The main instrument for the enumeration of unlabelled objects is the Burnside Lemma. It allows to express the number of unlabelled objects through the number of labelled objects which are fixed by the action of some group $G$ that defines the equivalence relation. In the paper \cite{Mednykh_Nedela} it has been shown that in our case we could restrict the action of $G$ to the set of periodic orientation-preserving homeomorphisms of the embedding surface. So we have to describe and enumerate the maps that are fixed under the action of the corresponding homeomorphisms. For describing these maps we will use the notion of so-called quotient maps on orbifolds \cite{Mednykh_Nedela}.

Consider as an example the representation of torus as a $2n$-polygon with the opposite sides pairwise identified. Assume that there is a possibility to draw some map on this polygon is such a way that it is fixed by the rotation of the polygon by the angle of $2\pi/L$. In this case we say that this map admits this automorphism of period $L$. The polygon's points with respect to the described action are divided into two subsets: the infinite subset of points in a general position and the finite subset of points in a special position (so-called branch points). Any point in general position lies on the orbit of period $L$, while a branch lies on an orbit with a period less than $L$. From the topological point of view any such torus $X$ could be viewed as a $L$-fold covering of a sphere $O$ with some branch points. In the general case this manifold $O$ with a finite number of branch points is called $g$-admissible orbifold, where $g$ is a genus of the covering manifold $X$. Any point of $O$ different from a branch point has the same number $L$ of preimages. Any branch point has the number of preimages less than $L$ and is characterized by so-called branch index $m$ equal to $L$ divided by the number of such preimages. Finally, the $L$-folded covering of a manifold $O$ by $X$ could be described by so-called signature
$$
(g,\mathfrak{g},L,[m_1,\ldots,m_r]),\qquad 1<m_1\leq\ldots\leq m_r,
$$
where $\mathfrak{g}$ is the genus of orbifold $O$, $r$ is the number of branch points, and $m_i$ are their branch indices.

Let us get back to the map $\M$ on the torus which is transformed onto itself under the rotations by the angle $2\pi/L$ and consider $1/L$-th part of this map. We could obtain the map $\M$ by gluing $L$ copies of such parts together. Each such part is a so-called quotient map $\mathfrak{M}$ drawn on the orbifold $O$. Suppose that one of the vertices $\mathfrak{x}$ of $\mathfrak{M}$ coincides with the branch point of $O$ having index $m_i$. Then this vertex corresponds to $L/m_i$ vertices of $\M$ with degrees $m_i$ times larger than the degree of $\mathfrak{x}$. Similarly, if a branch point falls into the center of some face of $\mathfrak{M}$, then this face is multiplied into $L/m_i$ faces of $\M$. Each of those faces has $m_i$ times more incident vertices and edges than the original face. 

A quotient map also differs from the ordinary map as it can contain dangling semi-edges besides the ordinary edges. They appear when the preimages of some branch point fall into the centers of edges of $\M$. It is necessary that this branch point has branch index $2$. Any ordinary edge of $\mathfrak{M}$ is multiplied into $L$ edges of $\M$, but a semi-edge corresponds only to $L/2$ edges of $\M$.

Note that several different homeomorphisms could correspond to the same orbifold. Mednykh and Nedela have proved that the number of them is given by the formula 
\begin{equation}
\label{eq:signat_number}
\Epi_0(\pi_1(O),\Z_L)=m^{2\mathfrak{g}}\phi_{2\mathfrak{g}}(L/m)\,E(m_1,\ldots,m_r),
\end{equation}
where $m=\lcm(m_1,\ldots,m_r)$ is the least common multiple of branch indices,
$$
E(m_1,\ldots,m_r)=\dfrac{1}{m}\sum\limits_{k=1}^m\prod\limits_{i=1}^r\Phi(k,m_i),\qquad \text{where}\qquad
\Phi(k,m_i)=\dfrac{\phi(m_i)}{\phi(n_i)}\,\mu(n_i),\qquad n_i=\dfrac{m_i}{\gcd(k,m_i)},
$$
$\phi(n)$ is the Euler function, $\mu(n)$ is the Moebius function. Besides that we should take into account the following additional constraints on the indices $m_i$ \cite[Th. 4.3]{Mednykh_Nedela}:
\begin{equation}
\label{eq:Riemann-Hurvitz}
2-2g=L\left(2-2\mathfrak{g}-\sum\limits_{i=1}^r\left(1-\dfrac{1}{m_i}\right)\right),
\end{equation}
\begin{equation}
\label{eq:H2}
m=\lcm(m_1,\ldots,m_r)\mid L;
\end{equation}
\begin{equation}
\label{eq:H3}
\lcm(m_1,\ldots,m_{i-1},m_{i+1},\ldots,m_r)=m\quad\forall\,\,i=1,\ldots,r.
\end{equation}

Summing up the above one could get from the Burnside's lemma the following formula for the number of unrooted maps on the surface of genus $g$ \cite[Th. 3.1]{Mednykh_Nedela}:
\begin{equation}
\label{eq:mednyh_main}
\tilde{\epsilon}(g) = \frac{1}{2n}\sum_{L|2n} \sum_{O \in Orb(S_g, Z_L)} \Epi_0(\pi_1(O),\Z_L) \epsilon_O(2n/L)
\end{equation}
Here $\epsilon_O(2n/L)$ is the number of rooted quotient maps with $2n/L$ semi-edges on an orbifold $O$ corresponding to maps with specified properties, in our case, $4$-regular one-face maps.

To use this formula we have to solve two problems. First, we have to describe all the orbifolds which satisfy the constraints (\ref{eq:Riemann-Hurvitz}), (\ref{eq:H2}) and (\ref{eq:H3}). Second, we have to understand how do quotient maps that correspond to $4$-regular one-face maps look like, and enumerate them.   

We note first of all that if the original map $\M$ on a surface $X$ was unicellular, then the map $\mathfrak{M}$ on an orbifold $O$ would also be unicellular. A face on an orbifold corresponds to a single face of the original map if and only if it contains a branch point of index $L$. In other words, we should take into account only those orbifolds that contain at least one branch point of index $L$. To be specific, let $m_r=L$. From the condition (\ref{eq:H2}) it follows that 
$$
m=\lcm(m_1,\ldots,m_{r-1},L)\mid L\qquad \Longrightarrow\qquad m=L.
$$

Now assume that all the vertices of one-faced map $\M$ have the same degree $d=4$. Note that all the other branch points of an orbifold $O$ should coincide either with the vertices or with the dangling ends of semi-edges of $\mathfrak{M}$. In the latter case the branch point should have the branch index equal to $2$. A branch point of an index $m_i$ could coincide with the vertex $\mathfrak{x}$ of $\mathfrak{M}$ only if $m_i$ divides $4$ and $\deg(\mathfrak{x})=4/m_i$.  

We noted above that $m=L$. Setting $i$ to be equal to $r$ in (\ref{eq:H3}) we get that
\begin{equation}
\label{eq:H3_1}
m=L=\lcm(m_1,\ldots,m_{r-1}).
\end{equation}
As all $m_i$, $i=1,\ldots,r-1$, divide $4$, it follows from (\ref{eq:H3_1}) that $L$ also divides $4$. The case $L=1$ corresponds to a trivial automorphism and all the labelled maps $\M$ on a surface $X$. The number of those is given by the sequence $\epsilon^{(4)}_g$ defined in (\ref{eq:num_4_maps}). It follows that for our case of 4-regular one-face maps the equality (\ref{eq:mednyh_main}) could be rewritten as
\begin{equation}
\label{eq:mednyh_specified}
\tilde{\epsilon}(g) = \frac{1}{4g-2} (\epsilon^{(4)}_g + f_2(g)+f_4(g)),
\end{equation}
where $f_2(g)$ and $f_4(g)$ are the summands corresponding to $L=2$ and $L=4$. 

Next we define ranges for variables in orbifold signatures and also count the numbers of corresponding epimorphisms. Consider first the case $L=2$. As it can be seen from the condition (\ref{eq:H3_1}), all $r$ branch points of the orbifold $O$ are in this case of the index $m_i=2$. Substitution of these values into the formula (\ref{eq:signat_number})  provides the number of epimorphisms equal to $m^{2\mathfrak{g}}=4^{\mathfrak{g}}$. The Riemann-Hurwitz equation (\ref{eq:Riemann-Hurvitz}) in the case $m_i=2$ becomes
$$
1-g=2-2\mathfrak{g}-r/2\qquad\Longleftrightarrow\qquad r=2g+2-4\mathfrak{g}.
$$
This equation and the inequality $r>0$ imply that $\mathfrak{g}$ varies form $0$ to $\lfloor g/2\rfloor$. 

Now let $L$ be equal to $4$. Condition (\ref{eq:H3_1}) means that among the branch points coinciding with the vertices of the factor-map $\mathfrak{M}$, there should be at least one point with the index $m_i=4$. We denote by $r_4\geq 1$ the number of branch points with the index $4$, and by $r_2$ number of branch points with the index $2$. Then the Riemann-Hurwitz equation (\ref{eq:Riemann-Hurvitz}) yields the following:
$$
2g=8\mathfrak{g}-6+3r_4+2r_2.
$$
From this equation it follows, first, that  $r_4$ is an even number. Second, $\mathfrak{g}$ varies from zero to some value corresponding to the minimum values of $r_4=2$ and $r_2=0$, that is, from zero to $\lfloor g/4\rfloor$. When $g$ and $\mathfrak{g}$ are fixed, the parameter $r_4$ varies from $r_4^{\min}=2$ to the value
$$
r_4^{\max}=\left\lfloor2(g+3-4\mathfrak{g})/3\right\rfloor
$$
that corresponds to $r_2=0$. Finally, from the formula (\ref{eq:signat_number}) we see that the number of epimorphisms for $L=4$ is equal to $2^{4\mathfrak{g}-1+r_4}$.

It follows from (\ref{eq:4-regular}) that any one-face $4$-regular map $\M$ on a surface $X$ with the genus $g$ has $2g-1$ vertices. Next we count the number of vertices of the corresponding factor-map $\mathfrak{M}$ on an orbifold $O$ with the genus $\mathfrak{g}$.  

Consider first a quotient map $\mathfrak{M}$ on an orbifold $O$ that corresponds to the case $L=2$ and has $r$ branch points of index $2$. Quotient map vertices that coincide with the branch points of the orbifold have degrees equal to 4. The remaining vertices of this map have degree 4. It is also convenient to think of dandling semi-edges of the quotient map $\mathfrak{M}$ as of normal edges of some other map $\tilde{\M}$ on a surface of genus $\mathfrak{g}$ with leafs at their ends. We denote by $s$ the number of leafs of the map $\tilde{\M}$ (that is, the number of danglind semi-edges of the quotient map $\mathfrak{M}$), by $l$ the number of vertices of degree $2$, and by $k$ the number of vertices of degree $4$. Next we obtain the relation between the parameters $s$, $l$ and $k$. 

Since the number of vertices of degree $2$ is equal to the number of branch points coinciding with the vertices of $\mathfrak{M}$, it follows that 
$$
l=r-1-s=2g+1-4\mathfrak{g}-s\qquad \Longrightarrow\qquad s=2g+1-4\mathfrak{g}-l.
$$
To calculate the number $k$ of vertices of degree $4$ it remains to subtract the number $l$ of vertices coinciding with the branch points from the number $2g-1$ of the vertices of the original map and then divide the result by two: each vertex of a factor-map $\mathfrak{M}$ on the orbifold that does not coincide with any of the branch points corresponds to $L = 2$ vertices of the original map $\M$. So,
$$
k=\dfrac{2g-1-l}{2}, \qquad \qquad l=2g-1-2k,\qquad \qquad s=2k-4\mathfrak{g}+2.
$$

We now turn to the description of the quotient-map $\mathfrak{M}$ corresponding to the map $\M$ in the case $L=4$. As noted above, in this case one branch point of index $4$ necessarily falls into the face. The remaining $r_4 - 1$ branch points of index $4$ necessary fall into the vertices of degree 1 of the quotient map $\mathfrak{M}$. Again, denote by $s$ the number of branch points coinciding with the dangling semi-edges of the quotient map $\mathfrak{M}$. As well as in the previous case, they could be viewed as vertices of a map $\tilde{\M}$ on the surface of genus $\mathfrak{g}$. Thus, the map $\tilde {\M} $ has a total of $\tilde{s} = s + r_4 - 1$ leafs.

Now we count the numbers $l$ and $k$ of vertices of degree $2$ and $4$ correspondingly. Vertices of degree $2$ are the ones that coincide with the orbifold's branch points of index $2$. Since $s$ out of $r$ such branching points fall into dangling semi-edges of the quotient map $\mathfrak{M}$, then the remaining $r_2 - s$ ones fall into vertices, so $l=r_2-s$. The number $k$ of vertices of degree 4 in this case is given by the formula
$$
k=\dfrac{2g-1-(r_4-1)-2(r_2-s)}{4}=\dfrac{2g-r_4-2r_2+2s}{4}=\dfrac{4\mathfrak{g}-3+r_4+s}{2}.
$$
Indeed, from the total amount $2g-1$ of vertices of the original map $\M$ we should subtract the amount $r_4-1$ vertices that turned into leafs on the quotient map, the doubled number $2(r_2-s)$ of the vertices that turned into vertices of degree 2, and then divide the result by 4: each vertex that does not coincide with the orbifold's branch points is multiplied into $L=4$ vertices on the surface of genus $g$. Consequently, 
$$
s=2k+3-4\mathfrak{g}-r_4,\qquad l=g-2k+r_4/2.
$$
From the first equality, in particular, it follows that $s$ is odd. In addition, the minimum value $s_{\min}$ of $s$ is $1$. As a result, $k_{\min}=2\mathfrak{g}-1+r_4/2$. The maximal value of $k$ is obtained when $s$ is also maximal. But $s=r_2-l$, that's why $s$ (and $k$) reach their maximal values when $l=0$, so $k_{\max}=\lfloor g/2+r_4/4\rfloor$. 

We turn, finally, to the calculation of the number of unlabelled 4-regular one-face maps. As shown above (see formula (\ref{eq:num_1_4_maps})), the number $\epsilon^{(1\div 4)}_\mathfrak{g}(k)$ of rooted $1\div 4$-valent one-face maps on the surface of genus $\mathfrak{g}$ that have $k$ vertices of degree $4$, $n$ edges and $s$ leafs, is calculated by the formula
$$
\epsilon^{(1\div 4)}_\mathfrak{g}(k)=\dfrac{2\,n!}{4^\mathfrak{g}\,\mathfrak{g}!\,s!\,(k-\mathfrak{g})!}=\dfrac{2(3k+1-2\mathfrak{g})!}{4^\mathfrak{g}\mathfrak{g}!(2k+2-4\mathfrak{g})!(k-\mathfrak{g})!}.
$$

In the case of $L = 2$ such maps have $n=\left[s+4k\right]/2=3k-2\mathfrak{g}+1$ edges, and among those we can distribute $l$ vertices of degree 2 in the number of ways equal to
$$
\BCCf{n}{l}=\BCf{n+l-1}{l}=\dfrac{(n+l-1)!}{(n-1)!\,l!}.
$$
However, the initial $(1\div 4)$-map had $2n$ ways to choose a root. In the new $(1\div 2\div 4)$-map a root could be chosen in $4g-2$ ways, as $4g-2$ semi-edges of the map $\tilde{\M}$ on an orbifold will become $8g-4$ semi-edges of the map $\M$ on the surface of genus $g$. Therefore, for the surface of genus $\mathfrak{g}$ and for $L=2$ we have 
$$
\epsilon^{(1\div 4)}_\mathfrak{g}(k)\dfrac{(n+l-1)!}{2\,n!\,l!}(4g-2)=\dfrac{(n+l-1)!}{4^\mathfrak{g}\,\mathfrak{g}!\,s!\,(k-\mathfrak{g})!\,l!}(4g-2)
$$
rooted $(1\div 2\div 4)$-maps $\tilde{\M}$. Multiplying this number by the number $4^{\mathfrak{g}}$ of epimorphisms and summing over $\mathfrak{g}$ and $k$, we get the expression  
$$
f_2(g)=(4g-2)\sum\limits_{\mathfrak{g},k}\dfrac{(n+l-1)!}{\mathfrak{g}!\,s!\,(k-\mathfrak{g})!\,l!}=
(4g-2)\sum\limits_{\mathfrak{g},k}\dfrac{(2g-2\mathfrak{g}+k-1)!}{\mathfrak{g}!\,(2k-4\mathfrak{g}+2)!\,(k-\mathfrak{g})!\,(2g-1-2k)!}.
$$
In this formula $\mathfrak{g}$ varies from $0$ to $\lfloor g/2\rfloor$, $k$ varies from $\max(0,2\mathfrak{g}-1)$ to $g-1$. Note that for $\mathfrak{g}=0$ the sum over $k$ could be rewritten as
$$
(4g-2)\sum\limits_{k=0}^{g-1}\dfrac{(2g+k-1)!}{(2k+2)!\,k!\,(2g-1-2k)!}=\dfrac{3}{2g+1}\BCf{4g-2}{2g}.
$$
For the other values of $\mathfrak{g}$ we failed to simplify the sum over $k$.

In the case of $L = 4$ the number $n$ of edges is $n=\left[\tilde{s}+4k\right] / 2= \left[s+r_4-1+4k\right] /2$.
Distributing the vertices of degree $2$ among these edges in the number of ways equal to
$$
\BCCf{n}{l}=\BCf{n+l-1}{l}=\dfrac{(n+l-1)!}{(n-1)!\,l!},
$$
we get some number of $(1\div 2\div 4)$-maps on an orbifold. However, original maps had $2n$ ways to distinguish a root semi-edge, but any new map should have $2g-1$ ways to do it: $2g-1$ semi-edges on the orbifold will become $8g-4$ semi-edges on the original surface. Therefore, we have to use double counting and multiply the fraction above by the fraction 
$$
\dfrac{2g-1}{2n}.
$$
The last thing that remains is to multiply the obtained expression by the binomial coefficient
$$
\BCf{\tilde{s}}{s}=\dfrac{\tilde{s}!}{(\tilde{s}-s)!\,s!},
$$
that counts the number of ways to choose among $\tilde{s}$ leafs those $s$ vertices that correspond to dangling semi-edges of the map $\mathfrak{M}$ on an orbifold. As a consequence, the number of quotient $(1\div 2\div 4)$-maps on an orbifold of genus $\mathfrak{g}$ by the formula
$$
\dfrac{(2g-1)\,(n+l-1)!}{4^\mathfrak{g}\,\mathfrak{g}!\,(k-\mathfrak{g})!\,l!\,s!\,(\tilde{s}-s)!}.
$$
Multiplying this fraction by the number $2^{4\mathfrak{g}-1+r_4}$ of epimorphisms, substituting into this formula the expressions for $n$, $l$, $\tilde{s}$ and $s$, and summing over $k$, $r_4$ and $\mathfrak{g}$, we obtain the function $f_4(g)$ in the form
$$
f_4(g)=
(2g-1)\sum\limits_{\mathfrak{g},r_4,k}\dfrac{2^{2\mathfrak{g}-1+r_4}\,(n+l-1)!}{\mathfrak{g}!\,(k-\mathfrak{g})!\,l!\,s!\,(\tilde{s}-s)!}=
$$
$$
=(2g-1)\sum\limits_{\mathfrak{g},r_4,k} \dfrac{2^{2\mathfrak{g}-1+r_4}\,(k-2\mathfrak{g}+g-r_4/2)!}{\mathfrak{g}!\,(k-\mathfrak{g})!\,(g-r_4/2-2k)!\,(2k+3-4\mathfrak{g}-r_4)!\,(r_4-1)!}.
$$

Substituting $\epsilon_g^{(4)}$, $f_2(g)$ and $f_4(g)$ into the formula (\ref{eq:mednyh_specified}), we finally get the following expression for the number $\tilde{\epsilon}^{(4)}(g)$ of $4$-regular unlabelled one-face maps of genus $g$:
$$
\tilde{\epsilon}^{(4)}(g)=\dfrac{(4g-3)!}{4^g\,g!\,(g-1)!}+\dfrac{3(4g-3)!}{2\,(2g+1)!\,(2g-2)!}+
$$
\begin{equation}
\label{eq:num_nonrooted_4}
+\sum\limits_{\mathfrak{g}=1}^{\lfloor g/2\rfloor}\sum\limits_{k=2\mathfrak{g}-1}^{g-1}\dfrac{(2g-2\mathfrak{g}+k-1)!}{2\,(2k-4\mathfrak{g}+2)!\,\mathfrak{g}!\,(k-\mathfrak{g})!\,(2g-1-2k)!}+
\end{equation}
$$
+\sum\limits_{\mathfrak{g}=0}^{\lfloor g/4\rfloor}\sum\limits_{\substack{r_4=2\\ 2 \mid r_4}}^{\lfloor2(g+3-4\mathfrak{g})/3\rfloor}
\sum\limits_{k=2\mathfrak{g}-1+r_4/2}^{\lfloor g/2+r_4/4\rfloor}
\dfrac{2^{2\mathfrak{g}-3+r_4}\,(k-2\mathfrak{g}+g-r_4/2)!}{\mathfrak{g}!\,(k-\mathfrak{g})!\,(g-r_4/2-2k)!\,(2k+3-4\mathfrak{g}-r_4)!\,(r_4-1)!}.
$$

\section{Conclusion}

In this work we obtained the explicit analytical expressions (\ref{eq:num_4_maps}) and (\ref{eq:num_nonrooted_4}) for counting $4$-regular one-face maps on a surface of an arbitrary genus $g$. The following table lists the numbers computed by those formulas.

\begin{table}[h!]
\footnotesize
\centering
\begin{tabular}{c|cc}
\midrule
$g$  &\phantom{000000000000000}Labelled\phantom{000000000000000}&\phantom{000000000000000}Unlabelled\phantom{000000000000000}\\ 
\midrule
1 & 1 & 1 \\
2 & 45 & 6 \\
3 & 9450 & 510 \\
4 & 4729725 & 169772 \\
5 & 4341887550 & 120644422 \\
6 & 6352181485650 & 144369379620 \\
7 & 13566444744352500 & 260893265836244 \\
8 & 39834473380605028125 & 663907896121296616 \\
9 & 153946961458244898693750 & 2263925904300525582790 \\
10 & 757572997336023146471943750 & 9968065754464730977513732 \\
11 & 4625189759553876588251163487500 & 55061782851836038471634743076 \\
12 & 34307345041490879593353005168531250 & 372905924364031740449809951518408 \\
13 & 303883906271359598859584503473567187500 & 3038839062713596039129776983675546524 \\
14 & 3168250194798584983481619521143486701562500 & 29335649951838749853328539549957507066456 \\
15 & 38405528861348447169764191835301345796340625000 & 331082145356452130774665205463914398071175024 \\
\midrule
\end{tabular}
\caption{Numbers of $4$-regular one-face maps by genus}
\end{table}

\newpage

\bibliography{Ref_book}

\end{document}